\documentclass[10pt,a4paper,oneside]{article}

 \RequirePackage{fixltx2e}            
 \RequirePackage[fleqn]{amsmath} 
 \RequirePackage{amsfonts}
 \RequirePackage{amssymb}
 \RequirePackage[all]{onlyamsmath}
 
 \RequirePackage[pdftex,hyperfootnotes=false,pdfpagelabels,pagebackref]{hyperref}
 \pdfcompresslevel=9 \pdfadjustspacing=1 \hypersetup{colorlinks=true,%
 urlcolor=blue,linkcolor=blue,citecolor=blue,linktocpage=false,pdfstartpage=1,pdfstartview=FitV,%
 breaklinks=true,pdfpagemode=UseNone,pageanchor=true,pdfpagemode=UseOutlines,plainpages=false,%
 bookmarksnumbered,bookmarksopen=true,bookmarksopenlevel=0,hypertexnames=true,pdfhighlight=/O,%
 pdftitle={The Lcm\{ 1,2,\dots,n \} as a Product of Sine Values Sampled Over the Points in Farey Sequences},%
 pdfauthor={Peter Luschny and Stefan Wehmeier},pdfsubject={},%
 pdfkeywords={least common multiple, sine product, Farey sequence, Gamma product, reflection formula}}
 
 \newcommand{\seqnum}[1]{\href{http://www.research.att.com/cgi-bin/access.cgi/as/~njas/sequences/eisA.cgi?Anum=#1}{{#1}}}
 
\begin{document}

 \begin{titlepage}
 \title{The Lcm$\left( 1,2,\dots,n \right)$ as a Product 
 of Sine Values Sampled Over the Points in Farey Sequences}
 \bigskip
 \author{Peter Luschny \and Stefan Wehmeier}
 \medskip
 \end{titlepage}
 \maketitle

 \begin{abstract}
 Some easily proven trigonometric formulae are given. 
 They lead to a shorter, alternate proof of a formula of G. Martin.
 \end{abstract}

 \section{A formula for the lcm of \{1,2,...,n\}} 

 Recently Greg Martin \cite{Martin} derived an interesting 
 formula for the least common multiple of 
 $\left\{1,2,\ldots,n \right\}$. In this paper, we give 
 an exposition of the proof in terms of the sine function.
 
 Let us first agree on some notation. 
 We write $\operatorname{LCM}(n)$ for 
 $\operatorname{lcm} \left\{1,2,\ldots,n \right\}$.
 ${F}(n)$ will denote the Farey sequence of order $n$, 
 that is the set of all reduced fractions in the interval 
 $\left[0,1\right]$ whose denominators are $n$ or less, 
 arranged in increasing order. We write $k \perp n$ 
 if $k$, $n$ are integers and relatively prime, following 
 Donald E. Knuth \cite{GKP}. We say an integer $n$ is 
 not a prime power if at least two different primes 
 divide $n$ and write $n \neq p^\alpha$ in this case; 
 we will always assume $\alpha$ an integer $ > 0$.
 Also note that we write $\prod^{2}{P}$ instead of
 $(\prod {P})^{2}$ to avoid big brackets around
 products.

 Martin first proves
\begin{equation}  
\prod_{\substack{0 < k < n \\ k \perp n }}
\frac{\Gamma^{2}\left(\frac{k}{n}\right)}{{2 \pi}} 
= \left\{ \begin{array}{rl}
1, &\mbox{if } n \neq p^{\alpha}; \\
\frac{1}{p}, &\mbox{if } n=p^{\alpha}.
\end{array} \right. \qquad (n \geq 2) 
\label{MartinTheorem1} \end{equation}
By multiplying for $n=2, 3, \ldots$, he derives
 \begin{equation}
 \operatorname{LCM}(n)=
 \prod_{\substack{r\in {F}(n) \\ 0<r<1 }}
 \frac{2\pi }{\Gamma^{2} \left(r\right)} .
\label{martin} \end{equation}

We observe that if $r$ is a member of ${F}(n)$ then
$1-r$ is also a member. By the reflection formula
of the Gamma function, we can trade in two Gamma 
evaluations for one sine evaluation
$\sin(\pi z) = \pi /(\Gamma(z)\Gamma(1-z)).$
This way Martin's formula becomes the simpler
\begin{equation}
\operatorname{LCM}(n)= \frac 12 \
\sideset{}{^{2}}\prod_{\substack{r\in {F}(n)
\\ 0<r\leq 1/2 }} 2 \sin\left(\pi r\right)
\qquad (n \geq 2)\ .
\label{lcmsin} \end{equation}

 In a discussion in the newsgroup de.sci.mathematik
 Jutta Gut \cite{Gut} observed that the left hand sides 
 of (\ref{lcmsin}) are equal for $n$ and $n-1$ if $n$ is 
 not a prime power; we may add the observation that 
 for $n>2$ the quotient of both equals $p$ if $n = p^{\alpha}.$ 

 Equivalently, we may apply the reflexion formula to Martin's theorem 
 (\ref{MartinTheorem1}) directly. This immediately gives
 \begin{equation}
 \prod_{\substack{0 < k < n \\ k \perp n }} 
 2 \sin \frac{\pi k}{n}
 = \left\{ \begin{array}{rl}
 1, &\mbox{if } n \neq p^{\alpha}; \\
 p, &\mbox{if } n=p^{\alpha}.
 \end{array} \right. 
 \qquad (n \geq 0) 
 \label{jutta} \end{equation}
 For $n > 2$ the range of the product can be reduced to 
 $1 \leq k \leq \left\lfloor n/2\right\rfloor$
 provided the product is raised to the square as
  $\sin \left({\pi k}/{n}\right)
  =\sin \left({\pi (n-k)}/{n}\right)$ .
 
\smallskip

 The interest we noted during these discussions motivates us
 to present an alternative: we first prove (\ref{jutta}),
 and then derive (\ref{martin}) from it.

 \section{Sines of roots of unity}
 For a short proof of (\ref{jutta}), we recall two well-known facts.
 
 \paragraph{Fact 1} is from elementary geometry: if the
 arc between two points on the unit circle has length
 $\theta$, then the length of the chord between them is
 $2 \sin(\theta/2)$. Applying this to the points $1$ and
 $\exp(\frac{2k\pi i}{n})$ where $-n/2 \leq k \leq n/2$ gives
 \begin{equation}
 \left| 1 - \exp\left(\frac{2k\pi i}{n}\right) \right|
 = 2 \left|\sin \frac{k\pi}{n}\right| .
 \label{sinChord} \end{equation}

 \paragraph{Fact 2} is about cyclotomic polynomials, 
 which we denote by $\Phi_n$. It can be found in many 
 standard texts on algebra; see, for example, \cite{Lang}, 
 p.280, Exercise 4. It says: if a prime $p$ divides $n$, then
 $\Phi_{np}(X) = \Phi_n(X^p)$; if $p$ does not divide $n$,
 then $\Phi_{np}(X) = \Phi_n(X^p) / \Phi_n(X)$. Plugging
 in $X=1$ and using induction gives for all $n > 2$
 \begin{equation} 
  \Phi_n(1) = \prod_{\substack{-n/2 \leq k \leq n/2 \\ k \perp n }}
 \left(1 - \exp\left(\frac{2k\pi i}{n}\right)\right)
 =  \left\{ \begin{array}{rl}
 1, &\mbox{if } n \neq p^{\alpha}; \\
 p, &\mbox{if } n=p^{\alpha}. \end{array} \right.
 \label{Phin1} \end{equation}

 \paragraph{Proof of (\ref{jutta}).} 
 The formula holds for $n=2$. 
 Taking absolute values in (\ref{Phin1}), the $k$-th and
 $(-k)$-th factor become equal.
 Using (\ref{sinChord}) gives (\ref{jutta}).

\bigskip

 Using $(X^n - 1)/(X-1)$ instead of a cyclotomic polynomial,
 the same method has been used recently on Planet Math
 \cite{PlanetMath} to give a concise proof of
\begin{equation}
\frac n{2^{n-1}}=\prod_{0<k<n}\sin \frac{\pi k}n \ .
\end{equation}
Clearly this can also be written as
\begin{equation}
n =\prod_{\substack{0<k<n \\ k \perp n }}
2\sin \frac{\pi k}{n}
\prod_{\substack{0<k<n \\ k \not{\perp} n }}
2\sin \frac{\pi k}{n} \ .
\label{sinsin} \end{equation}
From (\ref{sinsin}) and (\ref{jutta}) follows the 
counterpart of (\ref{jutta}).
\begin{equation}
\prod_{\substack{0<k<n \\ k \not{\perp} n }}
2\sin \frac{\pi k}{n} =
\left\{ \begin{array}{rl}
n/p, &\mbox{if } n = p^{\alpha}\,; \\
  n, &\mbox{if } n \neq p^{\alpha}. 
\end{array} \right. 
\qquad (n \geq 1) 
\label{dualsin} \end{equation}

These relations lead to a complementary form of Martin's
identities (\ref{MartinTheorem1}) and (\ref{martin}).
Let $\delta(n) = \{d: d \mid n \text{ and } 0<d<n \}$ denote 
the set of proper divisors of $n\geq 0$ and define 
$\operatorname{\overline{LCM}}(n)=\operatorname{lcm}(\delta(n))$
if $\delta(n)$ is not empty, otherwise $1$. Then 
\begin{equation}
\operatorname{\overline{LCM}}(n)=  
\prod_{\substack{ k \not{\perp} n 
\\ 0<k <n }} 2 \sin\frac{\pi k}{n}  
\quad (n \geq 1)\ .
\label{lcmtildesin} \end{equation}
Applying the reflexion formula of the $\Gamma$ function 
this can be rewritten as
\begin{equation}
\prod_{\substack{0 < k < n \\ k \not{\perp} n }}
\frac{{2\pi}}{\Gamma^{2}\left(\frac {k}{n} \right) }
= \operatorname{\overline{LCM}}(n)\, \qquad (n \geq 1) \; .
\label{lcmtilde} \end{equation}

The sequences $\operatorname{{LCM}}(n)$ and 
$\operatorname{\overline{LCM}}(n)$ are 
indexed in Sloane's \textit{Online Encyclopedia of 
Integer Sequence}s as 
\seqnum{A003418} and \seqnum{A048671} respectively. 

\section{Cosines of roots of unity}

 The same method works for cosines instead of sines. Our first,
 geometric fact then says that the chord between the points $-1$
 and $\exp({2\pi i k}/{n})$ has length $2 \cos( \pi k/ n)$,
 such that in the next step the cyclotomic polynomials must
 be evaluated at $-1$.

Let $\tilde{\epsilon}_{n}(k)=1 + \exp(2 \pi i k / n)$, then
 by induction it is easily proved from the recursion formulas that 
 \[
 \Phi_n(-1) = 
 \prod_{\substack{0 < k < n \\ k \perp n }}-\tilde{\epsilon}_{n}(k)
 = \left\{ \begin{array}{rl}
 p, &\mbox{if } n = 2p^{\alpha}; \\
 1, &\mbox{otherwise.} 
 \end{array} \right.
 \quad \qquad (n > 2)
 \]
 
 The zero factor $\Phi_2(-1)$ makes results boring,
 and we will avoid it in what follows. Since 
 $\tilde{\epsilon}_{n}(k)\tilde{\epsilon}_{n}(n-k)
 = \left(2 \cos \pi k / n \right)^2$ we get
 \begin{equation} \label{cosProd}
 \prod_{\substack{0 < k < n \\ k \perp n }}
 2 \left| \cos \frac{\pi k}{n} \right|
 = \left\{ \begin{array}{rl}
 p, &\mbox{if } n= 2 p^{\alpha}; \\
 1, &\mbox{otherwise}. \end{array} \right.
 \quad \qquad (n > 2)
 \end{equation}
  
 Multiplying (\ref{cosProd}) for all denominators
 below a given bound we obtain a similar result
 for Farey sequences as before:
 \[
 \sideset{}{^{2}}\prod_{\substack{r\in {F}(n)
 \\ 0<r < 1/2 }} 2 \cos\left(\pi r\right)  = \operatorname{LCM}(n/2) .
 \]

 \noindent In the case of cosines, too, the method can be
 applied easily to $(X^n-1)/(X-1)$ instead of a cyclotomic
 polynomial and gives
 \begin{equation}
 \prod_{1 \leq k \leq \left\lfloor n/2\right\rfloor}
 2 \cos \left(\frac{\pi k}{n}\right)
 =  \left\{ \begin{array}{rl}
 1, &\mbox{if } n  \mbox{ odd }; \\
 0, &\mbox{if } n \mbox{ even}. \end{array} \right.
 \end{equation}
 Again, the zeroes come from a zero factor $k=n/2$; the
 product of the other factors is $\frac{1}{2}n$ as can
 be seen by using $\frac{X^n - 1}{(X-1)(X+1)}
 = \sum_{m=0}^{(n/2) - 1} X^{2m}$ for even $n$.

 \section{The multiplication theorem revisited}
 Let us conclude with a formula that involves the gamma function again.
 We could construct it by converting our results on
 products of sines back to products of gammas (using the
 reflexion formula), but give another method here.

 The multiplication theorem of Gauss states
 \[
 \prod_{0\leq k\leq m-1}\Gamma \left( z+\frac km\right)
 = (2\pi )^{\frac 12\left( m-1\right)}
 m^{\left( 1/2-mz\right) }\Gamma \left( mz\right) \ .
 \]
 $ \text{The substitution } m \leftarrow \phi (n)+1,
 \ z \leftarrow \frac 1{\phi (n)+1}, \text{ leads to} $
 \begin{equation}
 \sqrt{\phi(n)+1}
 \prod_{0 < k \leq \phi(n)+1}\Gamma
 \left(\frac{k}{\phi(n)+1}\right)
 =(2\pi )^{\frac 12\phi (n)} \ .
 \label{multheo} \end{equation}
 On the other hand, if $\phi(n)$ denotes Euler's totient function and 
 dividing (\ref{martin}) for consecutive values $n$ and $n-1$, 
 we immediately see that if $n$ is not a prime power, then also
 \[
 \prod_{\substack{0 < k < n \\ k \perp n }}
 \Gamma \left(\frac{k}{n}\right)
 = (2\pi)^{\frac 12\phi (n)} \, .
 \]
 Abbreviating $N = \phi (n)+1$ and equating the left
 hand sides, we arrive at
 \[ 
 \prod_{\substack{0 < k < n \\ k \perp n }}
 \Gamma \left(\frac{k}{n}\right)
 =\sqrt{N}\prod_{ 0 < k < N}\Gamma \left(\frac{k}{N}\right)
 \quad (n \neq p^a) \   .
 \]
 It would be interesting to know whether there is some natural
 direct proof of this formula that does not use sines and 
 the reflexion formula.

\pdfbookmark[-1]{References}{References}

 \medskip
 \noindent 2010 Mathematics Subject Classification: 33B10 (11A05, 11B57).

 \bigskip
 \noindent E-mail address: \\
 \texttt{peter@luschny.de} \\
 \texttt{stefanw@math.upb.de}

 \end{document}